 \documentclass[12pt]{amsart}

\usepackage{graphics}
\usepackage{amssymb}

\newtheorem{theorem}{Theorem}[section]
\newtheorem{lemma}[theorem]{Lemma}
\newtheorem{corollary}[theorem]{Corollary}

\theoremstyle{definition}

\theoremstyle{remark}
\newtheorem{remark}[theorem]{Remark}

\numberwithin{equation}{section}

 \def\mcg{{\rm Mod}}

\begin{document}

 \title[Stable commutator length of a Dehn twist]
 {Stable commutator length of a Dehn twist}

 \author{Mustafa Korkmaz}
 \address{Department of Mathematics, Middle East Technical University,
 06531 Ankara, Turkey} \email{korkmaz@arf.math.metu.edu.tr}
 \subjclass{Primary 57M60, 20F38; Secondary 57N05, 57R17. }
 \date{\today}
 \keywords{Stable commutator length, Dehn twist, Mapping class group.}
 \begin{abstract}
 It is proved that the stable commutator length of a Dehn twist in 
the mapping class group is positive and
the tenth power of a Dehn twist about a nonseparating simple closed curve 
is a product of two commutators.  
As an application a new proof of the fact that the growth rate of
a Dehn twist is linear is given.
 \end{abstract}
 
 \maketitle

 \setcounter{secnumdepth}{1}
 \setcounter{section}{-1}

 \section{Introduction}

 The purpose of this paper is to prove that in the mapping class
 group of an orientable surface the stable commutator
 length of a Dehn twist about a simple closed curve not bounding a
 disc with punctures is positive for every $g$. We also give an
 upper bound for it. This gives an asymptotic estimate to
 Problem~$2.13$ $(B)(C)(D)$ in Kirby's problem book~\cite{k}.

 The upper bound for the stable commutator length of a Dehn
 twist is given based on known results. In the nonseparating case,
 however, we get a better upper bound by proving that the tenth
 power of such a Dehn twist is a product of two commutators, which
 is an interesting result itself.

 It was shown by Farb, Lubotzky and Minsky in \cite{fm} that
 the growth rate of a Dehn twist on an orientable surface of
 genus at least one is linear, answering a
 question of Ivanov (cf. Problem $2.16$ in \cite{k}).
 As an application of our main result we give a new proof
 of this fact by extending it to the genus zero case.

 The positivity of the stable commutator length of a Dehn twist
 about a separating simple closed curve was proved by Endo and
 Kotschick in~\cite{ek}. They concluded from this that the mapping
 class groups are not uniformly perfect and that the natural
 map from the second bounded cohomology to the ordinary
 cohomology of the mapping class group is not injective, which
 verified two conjectures of Morita. (cf. \cite{m},
 Conjectures~$6.19$ and~$6.21$.) In the proof of their main result,
 they use Seiberg-Witten theory.

 The main idea of our proof of the main result of this paper
 is to use the handlebody decomposition of a $4$-manifold
 admitting a Lefschetz fibrations. This was suggested
 to the author by Stipsicz for a signature computation
 in~\cite{korkmaz}.  In the computation we use the
 symplectic Parshin-Arakelov inequality of
 Li~\cite{li} to prove the nonexistence of certain Lefschetz fibrations.

 Donaldson~\cite{d} proved that every symplectic $4$-manifold
 admits a Lefschetz fibration after perhaps blowing up. Conversely,
 Gompf~\cite{gs} showed that the total space of every genus-$g$
 Lefschetz fibration admits a symplectic structure provided $g\geq 2$.
 This gives a combinatorial approach to symplectic $4$-manifolds
 through certain relations in mapping class groups. But
 understanding the relations in mapping class groups is not so easy.
 Usually, information in mapping class groups gives information about
 the corresponding $4$-manifolds. Examples of such applications are
 given in~\cite{ekks,ko,korkmaz}. The present paper, however,
 gives an application in the reverse direction.

 This paper has grown from a question of Andr\'as Stipsicz who asked
 to the author whether the mapping class groups were uniformly perfect.
 Author thanks him and Dieter Kotschick for their comments on
 the content of this paper, and the referee for his/her suggestions.

 \section{Preliminaries}

 For a compact orientable surface $S$ of genus $g$ with $p$ marked points 
 (to which we call punctures) and $q$ boundary components,
 we denote by $\mcg_{g,p}^q$ the
 mapping class group of $S$, the group of isotopy classes of
 orientation-preserving diffeomorphisms $S\to S$ which restrict
 to the identity on the boundary and preserve the set of punctures.
 The isotopies are also assumed to be the identity on the boundary
 and punctures. If $p$ and/or $q$ is zero, we omit it from the
 notation, so that, for example, $\mcg_g$ denotes $\mcg_{g,0}^0$.

 Let $S$ be an oriented surface.
 A simple closed curve $a$ on $S$ is called {\it trivial}
 if it bounds either a disc or a disc with one puncture.
 For every simple closed curve $a$ on $S$, there is a well known
 diffeomorphism called (right) Dehn twist about $a$, denoted by $t_a$,
 obtained by cutting the surface along $a$ and twisting one of the side to
 the right by $2\pi$ and gluing the two side back.
 A diffeomorphism and its isotopy class are denoted by the same symbol, and
 similarly for simple closed curves. For a Dehn twist $t_a$, we
 always assume that $a$ is nontrivial as the Dehn twist about a
 trivial simple closed curve is itself trivial in the mapping
 class group.

 We first state the next lemma which is elementary and will be
 used in the sequel. A proof of it can be found in~\cite{py}, Theorem~F, page xiii.   

 \begin{lemma}
 Let $\{a_n\}$ be a sequence of real numbers with nonnegative
 terms such that $a_{n+m}\leq
 a_n+a_m$ for every $n$ and $m$. Then the limit
 $\displaystyle \lim_{n\to\infty}
 \frac{a_n}{n}$ exists.
 \label{lem:sequence}
 \end{lemma}

 For a group $G$ let $[G,G]$ denote the commutator subgroup, the
 subgroup of $G$
 generated by all commutators $[a,b]=aba^{-1}b^{-1}$ for
 $a,b\in G$. For $x\in [G,G]$, we define
 the commutator length $c(x)$
 of $x$ to be the minimum number of factors needed to express
 $x$ as a product of commutators. Clearly,
 $c(x^{n+m})\leq c(x^{n}) + c(x^{m})$.  Therefore, we can define
 $$ ||x||= \lim_{n\to\infty}\frac{c(x^n)}{n},$$
 which is called the {\it stable commutator length of} $x$.

 Recall that for a group $G$, the first homology group $H_1(G)$ of
 $G$ with integer coefficients is isomorphic to the derived
 quotient group $G/[G,G]$.

 The next theorem and the corollary will be useful for us.

 \begin{theorem} [\cite{c}]
 Let $G$ be a group and let $u,v\in G$. Then $[u,v]^k$ can be written
 as a product of $E(\frac{k}{2})+1$ commutators,
 where $E(\frac{k}{2})$ denotes
 the integer part of $\frac{k}{2}$.
 \label{commute}
 \end{theorem}

 \begin{corollary} [\cite{b}] \label{cor:Culler-Bavard}
 Let $G$ be a group and let $u_1,v_1,u_2,v_2,\ldots u_r,v_r$ be elements of $G$.
 Then $([u_1,v_1][u_2,v_2]\ldots [u_r,v_r])^k$ is can be written
 as a product of $k(r-1)+E(\frac{k}{2}) +1$ commutators.
 \end{corollary}
 \begin{proof} The proof follows from
 $(uv)^k=(uvu^{-1})(u^{2}vu^{-2})\ldots
 (u^{k}vu^{-k})u^k$ and Theorem~\ref{commute}.
 \end{proof}

 \section{The main result}

 It is well known that the mapping class group $\mcg_g$ is
 perfect when $g\geq 3$ (cf.~\cite{p}). Hence, every element, in particular each
 Dehn twist, is a product of commutators.
 In the case of $g=2$, $H_1(\mcg_2)$ is isomorphic the
 cyclic group of order $10$ and is generated by the class of any Dehn
 twist about a nonseparating simple closed curve. A Dehn twist
 $t_a$ in $\mcg_2$ is not contained in the commutator subgroup, but $t_a^{10}$ is.
 Hence, we can talk about $||t_a^{10}||$ in this case.

Our main result is the following theorem.

 \begin{theorem}
 Let $S$ be a closed connected oriented surface of genus $g\geq 2$ and let
 $a$ be a nontrivial simple closed curve on $S$.
 Then $||t_a||\geq \frac{1}{18g-6}$ if $g\geq 3$ and
 $||t_a^{10}||\geq \frac{1}{3}$ if $g=2$.
 \label{th:||t_a||}
 \end{theorem}

 \begin{proof}
 Suppose first that $g\geq 3$. We assume the contrary that
 $||t_a||<\frac{1}{18g-6}$.
 Choose a rational number $r$ with $||t_a||< r<\frac{1}{18g-6}$.
 Then there exists an arbitrarily large positive integer $n$ such that
 $rn$ is an integer and $t_a^{n}$ can be written as a product
 of $rn$ commutators. This gives a relatively minimal
 genus-$g$ Lefschetz fibration over a closed orientable surface
 $\Sigma$ of genus $rn$ with the vanishing cycle $a$ repeated
 $n$ times as follows (We refer the reader to \cite{gs} for the details
 of the theory of
 Lefschetz fibrations).
 Consider $D^2\times S$, where $D^2$ is the $2$-disc.
 Attach $n$ $2$-handles to $\partial D^2\times S$ along the
 simple closed curve $a$ with $-1$ framing relative to
 the product framing. This gives a relatively minimal genus-$g$
 Lefschetz fibration
 $X_1\to D^2$ with monodromy $t_a^n$ along the boundary $\partial D^2$.
 Since $t_a^n$ is a product of $rn$ commutators, there is a surface bundle
 $X_2$ with fibers $S$ over an orientable surface of genus $rn$ with one boundary
 component such that the monodromy along the boundary is $t_a^n$.
 The boundary of $X_1$ and $X_2$ are genus-$g$ surface bundles over
 $S^1$ with monodromy $t_a^n$. Now glue $X_1$ and $X_2$ via an
 fiber preserving orientation reversing  diffeomorphism between boundaries
 to get a relatively minimal Lefschetz fibration $X\to \Sigma$ over a closed
 connected orientable surface $\Sigma$ with the generic fiber $S$.

 The Euler characteristic of $X$ is easily computed to be
 $$
 \chi(X)=4(g-1)(rn-1)+n=4grn-4rn-4g+4+n\ .
 $$
 Also, it follows from the construction that $b_1(X)\leq 2g+2rn$.

 We now give a lower bound for $b_2^-(X)$.
 For each $i=1,2,\ldots,n-1$, the cores of the $i$th and
 $(i+1)$st $2$-handles attached along $a$ give a sphere $S_i$ whose
 self intersection is $-2$. If $[S_i]$ denotes the homology class in 
 $H_2(X;\mathbb{R})$
 of  $S_i$, then $[S_i][S_{i+1}]=\pm 1$ and $[S_i][S_j]=0$ for $|i-j|\geq 2$,
 since the attaching regions of $2$-handles can be chosen to be disjoint.
 We can orient $S_i$ so that $[S_i][S_{i+1}]=-1$.
 It follows that the homology classes $[S_1],\ldots,[S_{n-1}]$ are linearly
 independent and form a basis for an $n-1$ dimensional subspace
 $V$ of $H_2(X;\mathbb{R})$. Hence, the matrix of the intersection form
 restricted to $V$ in the above basis is the matrix $-A$, where
 \[
 A=
 \begin{pmatrix}
  2 & 1 & 0 & 0 & \cdots & 0 & 0\\
  1 & 2 & 1 & 0 & \cdots & 0 & 0\\
  0 & 1 & 2 & 1 & \cdots & 0 & 0\\
  \cdot& \cdot & \cdot & \cdot & \cdots & \cdot & \cdot\\
  0 & 0 & 0 & 0 & \cdots & 2 & 1\\
  0 & 0 & 0 & 0 & \cdots & 1 & 2\\
 \end{pmatrix}.
 \]
 It is easy to check that the matrix $A$ is positive definite. Therefore, the
 restriction of the intersection form to $V$ is negative definite.
 It follows that
 \begin{eqnarray*}
 b_2^-(X)\geq n-1.
 \end{eqnarray*}

 An easy computation gives an upper bound for $b_2^+(X)$;
 \begin{eqnarray*}
 n-1 +b_2^+(X)
 &\leq& b_2(X)\\
 &=& \chi(X)+2b_1(X)-2\\
 &\leq & 4grn-4rn-4g+4+n+2(2g+2rn)-2\\
 &\leq & 4grn+n+2.
 \end{eqnarray*}
 Hence,
 \begin{eqnarray*}
 b_2^+(X)\leq 4grn+3.
 \end{eqnarray*}

 From these inequalities we obtain an upper bound for the signature of $X$;
 \[\sigma(X)\leq 4grn-n+4.\]

 On the other hand, following an argument of Kotschick~\cite{kot},
 Li proved in \cite{li} that
 $$ 2(g-1)(rn-1)\leq c_1^2(X).
 $$
 Hence, we obtain
 \begin{eqnarray*}
 2(g-1)(rn-1) & \leq & c_1^2(X)\\
 & = & 3\sigma (X)+2\chi (X) \\
 & \leq &  3(4grn-n+4)+2(4grn-4rn-4g+4+n)\\
 &  =   &  20grn-8rn-n+20-8g.
 \end{eqnarray*}
 As a result of this, we conclude that there exists arbitrarily big $n$ such that
 \begin{eqnarray*}
 0 & \leq &  [ (18g-6)r-1]n+18-6g.
 \end{eqnarray*}
 Since $(18g-6)r-1$ is negative, this is a contradiction.

 This proves the theorem for $g\geq 3$.

 If $S$ is a closed surface of genus two, then taking $n$ as a multiple of
 $10$ above finishes the proof.
 \end{proof}

 \begin{remark}
 I had originally proved that the stable commutator length of a
 Dehn twist in the above theorem are greater than or equal to
 $\frac{1}{20g-8}$, by using the inequality $c_1^2(X)\geq 0$
 proved in~\cite{s1}. The improvement was
 kindly suggested by Kotschick and Stipsicz.
  \end{remark}

 \begin{corollary}
 Let $S$ be a connected orientable surface of genus $g\geq 0$ with $p$
 punctures and $q$ boundary components such that $g+q\geq 2$.
 Let $a$ be a simple closed curve on $S$ not bounding a disc with punctures.
 Suppose that $t_a^k$ is in the commutator subgroup of
 $\mcg_{g,p}^q$. Then $||t_a^k||>0$.
 \label{cor:||t_a||}
 \end{corollary}

 \begin{proof}
 Let us glue a torus with one boundary component along each boundary
 component of $S$. By forgetting the punctures, we get a closed surface
 $R$ of genus $g+q$. The circle $a$ is now nontrivial on $R$.
 In this way we have a map $F$ from the mapping class
 group of $S$ to that of $R$. Clearly,
 $c(t_a^k)\geq c(F(t_a)^k)$. Since $F(t_a)$ is a Dehn twist about
 the nontrivial
 simple closed curve $a$ on $R$, the corollary follows from
 $||t_a^k||\geq ||F(t_a)^k||$ and  Theorem~\ref{th:||t_a||}.
 \end{proof}

 It was shown in~\cite{ko} that the commutator length of a Dehn
 twist is two. Theorem~\ref{th:||t_a||} and Theorem~\ref{commute}
 give the corollary.

 \begin{corollary}
 Let $S$ be a closed orientable surface of genus $g\geq 2$ and
 $a$ be a nontrivial simple closed curve on $S$.
 Then the element $t_a^k$ of $\mcg_g$ cannot be a commutator if
 $k > 9g-3$.
 \end{corollary}

 \begin{proof}
 Assume that $t_a^k$ is a commutator. Then $t_a^{kn}$ is a product
 of $E(\frac{n}{2})+1$ commutators, i.e. $c(t_a^{nk})\leq
 E(\frac{n}{2})+1$. Dividing both sides by $kn$ and taking the
 limit as $n$ tends to the infinity gives the desired contradiction
 $k \leq 9g-3$.
 \end{proof}

 \section{An upper bound for $||t_a||$.}

 In this section we give an upper bound for $||t_a||$
 for a simple closed curve $a$. In the case $a$ is nonseparating,
 we obtain a better upper bound by proving that the tenth power
 $t_a^{10}$ of a Dehn twist $t_a$ is a product of two commutators.

 The next lemma is well known.

 \begin{lemma} \label{braidrelation}
 Let $a$ and $b$ be two simple closed curves on an oriented surface $S$.
  \begin{itemize}
  \item[(a)] If $a$ is disjoint from $b$, then $t_a$ commutes with
  $t_b$.
  \item[(b)] If $a$ intersects $b$ transversely at only one point,
  then $t_at_bt_a=t_bt_at_b$.
 \end{itemize}
 \end{lemma}

 \begin{lemma}
 Let $S$ be a connected oriented surface and
 let $a,b,c$ and $d$ be four simple closed curves on $S$ such that
 there is an orientation preserving diffeomorphism of $S$
 mapping $a$ and $b$ to $d$ and $c$
 respectively. Then $t_at_b^{-1} t_ct_d^{-1}$ is a commutator.
 \label{lemma}
 \end{lemma}
 \begin{proof}
 Let $g$ be the isotopy class of a diffeomorphism mapping $a$ and
 $b$ to $d$ and $c$. Then
 \[
 t_at_b^{-1} t_ct_d^{-1}=t_at_b^{-1} t_{g(b)}t_{g(a)}^{-1}
 =t_at_b^{-1} gt_bt_a^{-1} g^{-1}
 =[t_at_b^{-1}, g] \ .
 \]
 \end{proof}

 \bigskip
 \begin{theorem}
 Let $S$ be a connected oriented surface of genus at least
 two and let $a$ be a nonseparating
 simple closed curve on $S$. Then $t_a^{10}$ can be written as a
 product of two commutators.
 \label{th:t_a^10}
 \end{theorem}

 \begin{proof}
 Since Dehn twists about two nonseparating simple closed curves
 are conjugate and a conjugate of a commutator is again a
 commutator, it suffices to prove the theorem for some nonseparating
 simple closed curve.

 Let $a_1,a_2,a_3$ be three nonseparating simple closed curves on
 $S$ such that $a_2$ intersects $a_1$ and $a_3$ transversely only
 once, $a_1$ is disjoint from $a_3$ and $a_1\cup a_3$ does not
 disconnect $S$. A regular neighborhood $a_1\cup a_2\cup a_3$
 is a torus with two boundary components, say $a_4$ and $a_5$,
 which are nonseparating on $S$.
 Clearly, $a_4$ and $a_5$ are disjoint from
 $a_1,a_2,a_3$, and from each other. Let us denote $t_{a_i}$ by $t_i$.
 It is well known that
 $t_4t_5=(t_1t_2t_3)^4$.  Using Lemma~\ref{braidrelation}, we obtain

 \begin{eqnarray*}
 t_4t_5
 &=&  (t_1t_2t_3)( t_1t_2t_3 )(t_1t_2t_3)( t_1t_2t_3)\\
 &=&  (t_1t_2t_1 )(t_3t_2t_3 )(t_1t_2t_1 )(t_3t_2t_3) \\
 &=&  (t_2t_1t_2)( t_2t_3t_2 )(t_2t_1t_2 )(t_2t_3t_2) \\
 \end{eqnarray*}
 Conjugating with $t_2^{-1}$ gives
 \begin{eqnarray*}
 t_4t_5
 &=&  t_1t_2t_2 t_3t_2t_2t_1t_2 t_2t_3t_2t_2 \\
 &=&  t_1 (t_2^2t_3t_2^{-2}) t_2^4
      t_1 t_2^{-1} (t_2^{3}t_3t_2^{-3})t_2^{-1} t_2^6\ .
 \end{eqnarray*}
 If we let $\alpha=t_{2}^2 (a_3)$ and $\beta=t_{2}^{3}(a_3)$,
 we get the equality
  \begin{eqnarray*}
    (t_4 t_{\alpha}^{-1} t_5t_1^{-1})
    = t_2^4(t_1 t_2^{-1} t_{\beta} t_2^{-1}) t_2^6\ .
  \end{eqnarray*}

 The curves $a_4$ and $a_5$ do not intersect $\alpha$ and $a_1$.
 Since the complements of $a_4\cup \alpha$ and $a_5\cup a_1$ are connected,
 there is an orientation preserving diffeomorphism
 taking $a_4$ and $\alpha$ to $a_1$ and
 $a_5$ respectively. The curve $a_2$ intersects $a_1$ and $\beta$
 transversely at one point. Hence, there is an orientation preserving
 diffeomorphism mapping $a_1$ and $a_2$ to $a_2$ and $\beta$ respectively.
 By Lemma~\ref{lemma}, each parenthesis is a commutator.
 Since the conjugate of a commutator is again
 a commutator, the proof follows.
 \end{proof}

 \begin{theorem}
 Let $S$ be a connected oriented surface of genus $g\geq 2$ and let
 $a$ be a simple closed curve on $S$. 
 \begin{itemize}
 \item[$(a)$] If $a$ is nonseparating, then $||t_a|| \leq \frac{3}{20}$ when $g\geq 3$ and
 $||t_a^{10}|| \leq \frac{3}{2}$ when $g=2$.
 
 \item[$(b)$] If $a$ is separating, then $||t_a|| \leq \frac{3}{4}$ when $g\geq 3$. 
 \end{itemize}
 \end{theorem}

 \begin{proof}
 If $a$ is nonseparating, by Theorem~\ref{th:t_a^10} and Corollary~\ref{cor:Culler-Bavard}, 
 the element $t_a^{10n}$ can be written as a product of
 $n+E(\frac{n}{2})  +1$ commutators.  The proof $(a)$ follows from this.

 If $a$ is separating, then $t_a^2$ is a product of two commutators~\cite{ko}.
It follows now from Corollary~\ref{cor:Culler-Bavard} that $||t_a||\leq \frac{3}{4}$.
 \end{proof}

 \begin{remark}
 In the case of an oriented surface of genus $2$, if $a$ is a separating simple 
 closed curve such that one of the components of the complement of $a$ has genus
 at two, then it is easy to conclude from the proof of
 Proposition~$10$ in~\cite{ko} that $t_a^2$ can be written a
 product of two commutators. It follows that $||t_a||\leq \frac{3}{4}$.

 If the genera of both components of the complement of $a$ are one,
 then $t_a^5$ is in the commutator subgroup of $\mcg_2$.  There are two nonseparating simple 
 closed curves $b,c$ on $S$ intersecting each other at one point such that 
 $t_a=(t_bt_c)^6$. It can be shown from this that $t_a^5$ is a product of $21$ commutators.
 It follows that $||t_a^5||\leq \frac{41}{2}$.
 \end{remark}

 \section{An application: Growth rate of Dehn twists}

 It is well known that the mapping class groups are finitely
 presented (cf.~\cite{w}).
 For a finite generating set $A$ of $\mcg_{g,p}^q$, let $d$ denote the corresponding
 word metric on $\mcg_{g,p}^q$. That is,  if $f,g\in \mcg_{g,p}^q$ then
 \[
 d(f,g)=\min \{n \,|\, g^{-1}f=x_1x_2\cdots x_n, x_i\in A\}.
 \]

 In \cite{fm}, Farb, Lubotzky and Minsky proved that if $g\geq 1$ then
 Dehn twists in the mapping class group $\mcg_{g,p}^q$  have linear growth
 rate. That is, the limit
 $$\displaystyle \lim_{n\to\infty} \frac{d(t_a^n,1)}{n}$$
 is positive.  Here, we give another proof of this fact and
 extend it to the genus zero case.
 Note that this limit does depend on the choice of the generating set,
 but the positivity of it does not.

 \begin{theorem}
 Let $S$ be a connected oriented surface of genus
 $g$ with $p$ punctures and $q$ boundary components.
 Suppose that $g+q\geq 2$. If $a$ is a simple closed curve on
 $S$ not bounding a disc with punctures, then for any finite generating set of
 the mapping class group, the limit
 $$\displaystyle \lim_{n\to\infty} \frac{d(t_a^n,1)}{n}$$ is positive,
 i.e. the growth rate of $d(t_a^n,1)$ is linear.
 \end{theorem}

 \begin{proof}
 Suppose first that $S$ is a closed surface so that $g\geq 2$.
 Let $A$ be a finite generating set for the mapping class group $\mcg_g$.
 In the case of $g=2$, we can choose $A$ so that it contains a set of generators for
 the commutator subgroup of $\mcg_2$.
 Let $k$ be a positive integer such that each element of $A$ contained
 in the commutator subgroup
 can be written as a product of $k$ commutators. Thus,
 $t_a^{10n}$ can be written as a product of $k\, d(t_a^{10n},1)$
 commutators for any positive integer $n$. Hence
$c(t_a^{10n})\leq kd(t_a^{10n},1)$.  It follows that
 $$
 \lim_{n\to\infty} \frac{d(t_a^n,1)}{n}\geq
 \frac{1}{10k}\, ||t_a^{10}||.$$
 Hence, it is positive by Theorem~\ref{th:||t_a||}.

 In the general case, let us glue a torus with one boundary to $S$
 along each boundary component of $S$ and forget the punctures.
 Let $A$ be a finite generating set for $\mcg_{g,p}^q$. As in the proof of
 Corollary~\ref{cor:||t_a||}, this gives a homomorphism
 $F:\mcg_{g,p}^q\to \mcg_{g+q}$. Extend $F(A)$ to a finite generating set $B$ for
 $\mcg_{g+q}$. Clearly, we have $d(t_a^n,1)\geq d(F(t_a)^n,1)$
 with respect to the generating sets $A$ and $B$.
 Since $F(t_a)$ is a Dehn twist about a nontrivial simple closed
 curve on a closed surface, the proof follows from the closed case.
 \end{proof}

 \section{A question}
 We end with a question, which arises from the topology of
 the Stein fillings of a contact $3$-manifold: Let $S$ be an
 oriented surface with one boundary component and let
 $a_1,a_2,\ldots$ be a sequence of nonseparating simple
 closed curves on $S$.  Does there exist a positive integer
 $N$ such that
 $$c(t_{a_1}t_{a_2}\cdots t_{a_n}) \geq Nn.$$
More generally, does the limit
 $$\lim_{n\to\infty} \frac{c(t_{a_1}t_{a_2}\cdots t_{a_n})}{n}$$
 exist? If it does, is it positive (for any choice of $a_i$)?
 Note that if all $a_n$ are equal, then we proved that there
 exists such an $N$.

 \bibliographystyle{amsplain}

\begin{thebibliography}{10}

 \bibitem{b}
 C.~Bavard,
 {\em Longueur stable des commutateurs},
 L'Enseignement Math\'{e}matique {\bf37} (1991), 109-150.

 \bibitem{c}
 M.~Culler,
 {\em Using surfaces to solve equations in free groups},
 Topology {\bf20} (1981), 133-145.

 \bibitem{d}
 S.~Donaldson,
 {\em Lefschetz fibrations in symplectic geometry},
 Proc. Internat. Cong. Math. (Berlin, 1998), Vol II, Doc. Math. Extra Volume
 ICM II (1998), 309-314.

 \bibitem{ekks}
 H.~Endo, M.~Korkmaz, D.~Kotschick, B.~Ozbagci, A.~I.~Stipsicz,
 {\em Commutators, Lefschetz fibrations and the signatures of surface bundles},
 Topology {\bf 41} (2002), 961-977.

 \bibitem{ek}
 H.~Endo, D.~Kotschick,
 {\em Bounded cohomology and non-uniform perfection of mapping class groups},
 Invent. Math. 144 (2001), 169--175.

 \bibitem{fm}
 B.~Farb, A.~Lubotzky, Y.~N.~Minsky,
 {\em Rank one phenomena of mapping class groups},
 Duke Math Journal, Vol. 106, (2001), 581-597.

 \bibitem{gs}
 R.~E.~Gompf and A.~I.~Stipsicz,
 {\em  4-manifolds and Kirby calculus},
 AMS Graduate Studies in Mathematics {\bf20} 1999.

 \bibitem{k}
 R.~Kirby, {\em Problems in low-dimensional topology},
 in Geometric Topology (W.~Kazez ed.) AMS/IP Stud. Adv. Math. vol 2.2,
 American Math.~Society, Providence 1997.

 \bibitem{korkmaz}
 M.~Korkmaz,
 {\em Noncomplex smooth 4-manifolds with Lefschetz fibrations},
 Internat. Math. Res. Notices  {\bf 2001}, 115-128.

 \bibitem{ko}
 M.~Korkmaz, B.~Ozbagci,
 {\em Minimal number of singular fibers in a Lefschetz fibration},
 Proc. Amer. Math. Soc. {\bf 129} (2001) 1545-1549.

 \bibitem{kot}
 D.~Kotschick,
 {\em Signatures, monopoles and mapping class groups},
 Math. Res. Lett. {\bf5} (1998), 227-234.

 \bibitem{li}
 T.~J.~Li,
 {\em Symplectic Parshin-Arakelov inequality},
 Internat. Math. Res. Notices {\bf2000}, 941--954.

 \bibitem{mm}
 S.~Matsumoto, S.~Morita,
 {\em Bounded cohomology of certain groups of homeomorphisms},
 Proc. Amer. Math. Soc. {\bf94} (1985), 539-544.

 \bibitem{m}
 S.~Morita,
 {\em Structure of the mapping class groups of surfaces: a survey
 and a prospect},
 Geometry and Topology Monographs Vol. 2, Proc. of Kirbyfest 1999,
 488-505.

 \bibitem{py}
M.~Pollicott, M.~Yuri, 
{\em Dynamical systems and ergodic theory}, 
London Mathematical Society Student Texts, 40. 
Cambridge University Press, Cambridge, 1998. 

 \bibitem{p}
J.~Powell
{\em Two theorems on the mapping class group of a surface}, 
Proc. Amer. Math. Soc. {\bf 68} (1978), 347--350.

 \bibitem{s1}
 A.~I.~Stipsicz,
 {\em Chern numbers of certain Lefschetz fibrations},
 Proc. Amer. Math. Soc. {\bf128} (2000), 1845-1851.

 \bibitem{w}
 B.~Wajnryb,
 {\em An elementary approach to the mapping class group of a surface},
 Geometry and Topology {\bf3} (1999), 405-466.

 \end{thebibliography}

 \end{document}